\documentclass[10pt,twoside]{amsart}
\def\YEAR{\year}\newcount\VOL\VOL=\YEAR\advance\VOL by-1995
\def\firstpage{1}\def\lastpage{1000}
\def\received{}\def\revised{}
\def\communicated{}

\makeatletter
\def\magnification{\afterassignment\m@g\count@}
\def\m@g{\mag=\count@\hsize6.5truein\vsize8.9truein\dimen\footins8truein}
\makeatother

\font\eightrm=cmr8
\font\caps=cmcsc10                    % Theorem, Lemma etc
\font\Caps=cmcsc10 scaled \magstep1   % Title

\usepackage{amsmath}
\usepackage{amssymb}
\usepackage{amscd}
\usepackage{color} 
\usepackage{graphicx}
\input xy 
\xyoption{all}

\makeatletter
\@specialpagefalse\headheight=8.5pt
\def\DocMath{}
\renewcommand{\@evenhead}{%
    \ifnum\thepage>\lastpage\rlap{\thepage}\hfill%
    \else\rlap{\thepage}\slshape\leftmark\hfill{\caps\SAuthor}\hfill\fi}%
\renewcommand{\@oddhead}{%
    \ifnum\thepage=\firstpage{\DocMath\hfill\llap{\thepage}}%
    \else{\slshape\rightmark}\hfill{\caps\STitle}\hfill\llap{\thepage}\fi}%
\makeatother

\def\TSkip{\bigskip}
\newbox\TheTitle{\obeylines\gdef\GetTitle #1
\ShortTitle  #2
\SubTitle    #3
\Author      #4
\ShortAuthor #5
\EndTitle
{\setbox\TheTitle=\vbox{\baselineskip=20pt\let\par=\cr\obeylines%
\halign{\centerline{\Caps##}\cr\noalign{\medskip}\cr#1\cr}}%
	\copy\TheTitle\TSkip\TSkip%
\def\next{#2}\ifx\next\empty\gdef\STitle{#1}\else\gdef\STitle{#2}\fi%
\def\next{#3}\ifx\next\empty%
    \else\setbox\TheTitle=\vbox{\baselineskip=20pt\let\par=\cr\obeylines%
    \halign{\centerline{\caps##} #3\cr}}\copy\TheTitle\TSkip\TSkip\fi%
\centerline{\caps #4}\TSkip\TSkip%
\def\next{#5}\ifx\next\empty\gdef\SAuthor{#4}\else\gdef\SAuthor{#5}\fi%
\ifx\received\empty\relax
    \else\centerline{\eightrm Received: \received}\fi%
\ifx\revised\empty\TSkip%
    \else\centerline{\eightrm Revised: \revised}\TSkip\fi%
\ifx\communicated\empty\relax
    \else\centerline{\eightrm Communicated by \communicated}\fi\TSkip\TSkip%
\catcode'015=5}}\def\Title{\obeylines\GetTitle}
\def\Abstract{\begingroup\narrower
    \parskip=\medskipamount\parindent=0pt{\caps Abstract. }}
\def\EndAbstract{\par\endgroup\TSkip}

\long\def\MSC#1\EndMSC{\def\arg{#1}\ifx\arg\empty\relax\else
     {\par\narrower\noindent%
     2000 Mathematics Subject Classification: #1\par}\fi}

\long\def\KEY#1\EndKEY{\def\arg{#1}\ifx\arg\empty\relax\else
	{\par\narrower\noindent Keywords and Phrases: #1\par}\fi\TSkip}

\newbox\TheAdd\def\Addresses{\vfill\copy\TheAdd\vfill
    \ifodd\number\lastpage\vfill\eject\phantom{.}\vfill\eject\fi}
{\obeylines\gdef\GetAddress #1
\Address #2 
\Address #3
\Address #4
\EndAddress
{\def\xs{5.0truecm}\parindent=0pt
\setbox0=\vtop{{\obeylines\hsize=\xs#1\par}}\def\next{#2}
\ifx\next\empty % 1 address
     \setbox\TheAdd=\hbox to\hsize{\hfill\copy0\hfill}
\else\setbox1=\vtop{{\obeylines\hsize=\xs#2\par}}\def\next{#3}
\ifx\next\empty % 2 addresses
     \setbox\TheAdd=\hbox to\hsize{\hfill\copy0\hfill\copy1\hfill}
\else\setbox2=\vtop{{\obeylines\hsize=\xs#3\par}}\def\next{#4}
\ifx\next\empty\ % 3 addresses
     \setbox\TheAdd=\vtop{\hbox to\hsize{\hfill\copy0\hfill\copy1\hfill}
                \vskip20pt\hbox to\hsize{\hfill\copy2\hfill}}
\else\setbox3=\vtop{{\obeylines\hsize=\xs#4\par}}
     \setbox\TheAdd=\vtop{\hbox to\hsize{\hfill\copy0\hfill\copy1\hfill}
	        \vskip20pt\hbox to\hsize{\hfill\copy2\hfill\copy3\hfill}}
\fi\fi\fi\catcode'015=5}}\gdef\Address{\obeylines\GetAddress}

\def\b#1{\bf #1}                   % Bold letters
\def\c#1{\mathcal #1}        % Math calligraphy
\def\op#1{\expandafter\def\csname #1\endcsname{\mathop{\mathrm{#1}}\nolimits}} % Define operators

\theoremstyle{plain}
    \newtheorem{theorem}{Theorem}[section]
    \newtheorem{proposition}[theorem]{Proposition}
    
    \newtheorem{corollary}[theorem]{Corollary}
    
\theoremstyle{definition}

\op{Gal}  \op{Hom} \op{End}  \op{Ext}   \op{Tor}  \op{Sym}  \op{Fil}
\op{Re}
\op{sing}   \op{can}  \op{cris}     \op{rig}  \op{an}    \op{st}   \op{univ}  \op{dR} \op{syn}
\op{Spec}   \op{Spf}   \op{Spm}
\op{pol}  \op{Im} \op{id}
     
\newcommand{\C}{\mathbb{C}}

\newcommand{\G}{\mathbb{G}}
\newcommand{\Z}{\mathbb{Z}}
\newcommand{\Q}{\mathbb{Q}}

\newcommand{\Zp}{{\mathbb{Z}_p}}

\newcommand{\zetapn}{\zeta_{p^n}}
\newcommand{\zetapnpi}{\zeta_{p^{n+1}}}

\newcommand{\zetapmpi}{\zeta_{p^{m+1}}}

\newcommand{\m}{\mathfrak{m}}

\begin{document}
%%%%% ------------- fill in your data below this line  -------------------
%%%%%    The following lines \Title ... \EndAddress must ALL be present
%%%%%    and in the given order.
\Title An elementary proof 
of  the Mazur-Tate-Teitelbaum conjecture
 for elliptic curves
%%%%%    Put here the title. Line breaks will be recognized. 
\ShortTitle   the Mazur-Tate-Teitelbaum conjecture
%%%%%    Running title for odd numbered pages, ONE line, please. 
%%%%%    If none is given, \Title will be used instead.          
\SubTitle   
Dedecated to Professor John Coates 
on the occasion of his sixtieth birthday
%%%%%    A possible subtitle goes here.
\Author Shinichi Kobayashi 
\footnote[1]{Supported by JSPS Postdoctoral Fellowships for Research Abroad.}
%%%%%    Put here name(s) of authors. Line breaks will be recognized.  
\ShortAuthor  Shinichi Kobayashi
%%%%%%   Running title for even numbered pages, ONE line, please. 
%%%%%%   If none is given, \Author will be used instead.          
\EndTitle  
\Abstract 
%%%%%    Put here the abstract of your manuscript.
%%%%%    Avoid macros and complicated TeX expressions, as this is
%%%%%    automatically translated and posted as an html file.
We give an elementary proof of  the Mazur-Tate-Teitelbaum conjecture for elliptic curves 
by using Kato's element. 
\EndAbstract
\MSC 
%%%%%    2000 Mathematics Subject Classification: 
11F85, 11G05, 11G07, 11G40, 11S40.
\EndMSC
\KEY
%%%%%    Keywords and Phrases:     
 elliptic curves, $p$-adic $L$-functions,  Iwasawa theory, 
the Mazur-Tate-Teitelbaum conjecture, 
exceptional zeros, Kato's element. 
\EndKEY
%%%%%    All 4 \Address lines below must be present. To center the last
%%%%%    entry, no empty lines must be between the following \Address
%%%%%    and \EndAddress lines.
\Address  
%%%%%    Address of first Author here
Shinichi KOBAYASHI 
Graduate School of Mathematics 
Nagoya University
Furo-cho Chikusa-ku
Nagoya 464-8602
Japan
shinichi@math.nagoya-u.ac.jp
\Address
%%%%%    Address of second Author here etc.
\Address
\Address
\EndAddress
%%
%%       Make sure the last tex command in your manuscript
%%       before the first \end{document} is the command  \Addresses
%%
%%---------------------Here the prologue ends---------------------------------
%%--------------------Here the manuscript starts------------------------------

\section{Introduction} 
The $p$-adic $L$-function $L_p(E,s)$ of an elliptic curve $E$ 
defined over $\Q$ has an extra zero at $s=1$  coming  from the interpolation factor at $p$ 
if $E$ has split multiplicative reduction at the prime $p$.  
The Mazur-Tate-Teitelbaum conjecture (now a theorem of Greenberg-Stevens) 
describes the first derivative 
of $L_p(E,s)$ as 
$$\frac{d}{ds}L_p(E,s)\; \vert_ {\, s=1}= \frac{\log_p (q_E)}{\mathrm{ord}_p (q_E)}\;\frac{L(E,1)}{\Omega^+_E} $$
where $q_E$ is the Tate period of $E$ coming from the $p$-adic uniformization of $E$ at $p$,   $\log_p$ is the Iwasawa $p$-adic logarithm, $\Omega^+_E$ is the real period of $E$ and  
$L(E, 1)$ is the special value of the complex Hasse-Weil $L$-function at $s=1$. 

Known proofs of this conjecture are classified into two kinds. 
One is, as Greenberg-Stevens \cite{GS} did first, 
a proof using a global theory like Hida's   universal ordinary deformation. 
The other is, as Kato-Kurihara-Tsuji \cite{KKT} or Colmez \cite{Co} did, 
a proof based on local theory (except using Kato's element). 
Each kind of proof has its own importance but the latter type of proof makes it clear that 
the substantial facts behind this conjecture are of local nature.  The $p$-adic $L$-function is the image of  Kato's element 
via a purely local morphism, the so called Coleman map or Perrin-Riou map. 
%the properties of the $p$-adic $L$-function are divided into  two kinds, 
%the properties coming from that of global Kato's element or 
% the properties coming from that of the local Coleman map. 
The extra zero phenomena discovered by  Mazur-Tate-Teitelbaum is, 
in fact, a property of the local Coleman map. 

In this paper,  we prove a derivative formula 
(Theorem \ref{dcol}) of the Coleman map for elliptic curves 
by purely local and elementary method and we apply this formula 
to Kato's element to show the conjecture of Mazur-Tate-Teitelbaum. 
Of course, our proof is just a special and the simplest case of that in 
Kato-Kurihara-Tsuji \cite{KKT} or Colmez \cite{Co}
(they proved 
the formula not only for elliptic curves but for higher weight 
modular forms) but 
I believe that it is still worthwhile to write it down for the following reason. 
First, the important paper Kato-Kurihara-Tsuji \cite{KKT} 
has not yet been published. Second, since we restrict ourselves to the case of 
elliptic curves, the proof is much simpler and elementary 
(of course, such a simple proof would be also known to specialists. 
In fact, Masato Kurihara informed me that Kato, Kurihara and Tsuji  have 
two simple proofs and one is similar to ours). 
I hope that this paper would help those who are interested in the understanding of 
this interesting problem. 
\ \\

\noindent {\b Acknowledgement}: 
I would like to wish Professor John Coates a happy sixtieth birthday, and 
to thank him for his contribution to mathematics, especially to Iwasawa theory. 
It is my great pleasure to dedicate this article to him on this occasion. 
 
 This paper was written during the author's visit at  the university of 
Paris 6. He would like to thank P. Colmez and L. Merel for the accommodation. 
He also would like to thank K. Bannai and N. Otsubo for  discussion. 
Finally, he is grateful to the referee for his careful reading of the manuscript.

\section{A structure of the group of local units in $k_\infty/\Q_p$.}

Let $k_\infty/\Q_p$ be the (local) cyclotomic $\Z_p$-extension 
in $\Q_p(\zeta_{p^\infty}):=\cup_{n=0}^\infty \Q_p(\zetapn)$
with Galois group $\Gamma$ and 
let $k_n/\Q_p$ be its $n$-th layer with 
Galois group $\Gamma_n$. 
We identify the Galois group $\mathrm{Gal}(\Q_p(\zeta_{p^\infty})/\Q_p)$ 
with $\Z_p^\times$ by the cyclotomic character $\kappa$. 
Then $\Gamma$ is identified with $1+p\Z_p$ and  
the torsion subgroup $\Delta$ of $\mathrm{Gal}(\Q_p(\zeta_{p^\infty})/\Q_p)$ 
 is regarded as $\mu_{p-1} \subset \Z_p^\times$. 

Let $U^1_n$ be the subgroup of $\c{O}_{k_n}^\times$ 
consisting of  the elements which are congruent to $1$ modulo the maximal ideal 
$\mathfrak{m}_n$ of $\c{O}_{k_n}$.  

Following the Appendix of  Rubin \cite{R} or \cite{K}, for a fixed generator $(\zetapn)_{n\in \mathbb{N}}$ of $\Zp(1)$,  
we construct a certain canonical  system of local points $(d_n)_n \in \varprojlim_n  U^1_n$ 
% following \cite{R} or  \cite{K}
and we determine the Galois module structure of $U^1_n$ using these points.

We let 
$$\ell(X)=\log(1+X)+\sum_{k=0}^\infty\; \sum_{\delta \in \Delta}\; 
\frac{(X+1)^{p^k \delta}-1}{p^k}.$$ 
This power series is convergent 
in $\Q_p[[X]]$ due to the summation $\sum_{\delta \in \Delta}$. 
%the $i$-th coefficient is equal to 
%$\frac{1}{p^k}\sum_k \sum_{\delta \in \Delta} \binom{p^k\delta }{i}$
%$\sum_{\delta}\; \frac{(X+1)^{p^k \delta}-1}{p^k} \equiv 0 \mod p^k \Zp[X].$ 
It is straightforward to see that 
$$\ell'(X) \in 1+X\Zp[[X]], \quad \ell(0)=0, \quad (\varphi -p) \circ \ell(X) \in p \Zp[[X]]$$ where $\varphi$ is 
the Frobenius operator such that $(\varphi \circ \ell)(X)=\ell((X+1)^p-1)$. 
Hence by Honda's theory, there is a formal group $\mathcal{F}$ over $\Z_p$ whose logarithm is 
given by $\ell$, and $\iota(X)=\exp \circ\, \ell \,(X)-1 \in \Zp[[X]]$ gives an isomorphism 
of formal groups 
$\mathcal{F} \cong \widehat{\G}_m$ over $\Z_p$. (See for example, Section 8 of \cite{K}.) 
Take an element $\varepsilon$ of $p\Zp$ such that $\ell (\varepsilon)=p$ and 
we define 
$$c_n\,:=\,\iota((\zetapnpi-1)\,[+]_\mathcal{F}\,\varepsilon).$$ 
Since this element is fixed by the group $\Delta$, 
this is an element of $\widehat{ \mathbb{G}}_m(\mathfrak{m}_n)$. 
Then by construction, $d_n=1+c_n \in U^1_n$ satisfies the relation 
$$\log_p (d_n) =\ell(\varepsilon)+\ell(\zetapnpi-1)=p+\sum_{k=0}^{n}\; \sum_{\delta \in \Delta}\; 
\frac{{\zeta_{p^{n+1-k}}}^{\delta}-1}{p^k}. $$ 

\begin{proposition}\label{norm}
i) $(d_n)_n$ is a norm compatible system and $d_0=1$. \\
ii) Let $u$ be a generator  of $U^1_0$. Then 
 as $\Z_p[\Gamma_n]$-module, 
$d_n$ and $u$ generate $U^1_n$, 
 and 
%In particular, we have 
%$$\mathrm{N}_{k_n/\Q_p} 
%U^1_n \,=\, ( U^1_0)^{p^n}\,=\,1+p^{n+1}\Zp,$$ and 
$d_n$ generates    $(U^1_n)^{\mathrm{N} =1}$ 
where $\mathrm{N}$ is the absolute norm from $k_n$ to $\Q_p$.
\end{proposition}
\begin{proof}  
Since $\zeta_p-1$ is not contained in $\mathfrak{m}_n$, 
the group $\widehat{\G}_m(\mathfrak{m}_n)$ does not contain $p$-power torsion points.  
Therefore to see i), it suffices to show the trace compatibility of $(\log_p (d_n))_n$, and 
this is done by direct calculations. 
For ii), we show that 
$(\iota^{-1}(c_n)^\sigma)_{\sigma \in \Gamma_n }$ and $\varepsilon$ generate $\mathcal{F}(\mathfrak{m}_n)$ as $\Z_p$-module by induction for $n$. 
The proof is the same as that of Proposition 8.11 of \cite{K} but 
we rewrite it for the ease of the reader. 
The case $n=0$ is clear.  
For arbitrary $n$, we show that 
$\ell (\mathfrak{m}_n)\subset \mathfrak{m}_n+k_{n-1}$ and 
$$\mathcal{F}(\mathfrak{m}_n)/\mathcal{F}(\mathfrak{m}_{n-1}) \cong \ell  
(\mathfrak{m}_n)
/\ell (\mathfrak{m}_{n-1})  \cong
 \m_n/\m_{n-1}.$$  
The first isomorphism is induced by the logarithm  $\ell$  and the 
last isomorphism is  by $(\m_n+k_{n-1})/k_{n-1} \cong \m_n/\m_{n-1}$. 
As a set, $\mathcal{F}(\mathfrak{m}_n)$ is the maximal ideal $\mathfrak{m}_n$, and  
we write $x \in \mathcal{F}(\m_n)$ in the form 
$x=\sum_{\delta \in \Delta}\sum_{i} a_i \, \zetapnpi^{i\delta}$, $a_i \in \Zp$. 
Then for  
$y=\sum_{\delta \in \Delta} \sum_{i} a_i \, \zetapn^{i \delta} \in \m_{n-1}$, 
we see that $x^{p} \equiv y \mod{p\c{O}_{k_n}}$ . 
Therefore  for $k \geq 1$, we have 
$$\sum_{\delta \in \Delta} \frac{(x+1)^{p^{k}\delta}-1}{p^k} \equiv \sum_{\delta \in \Delta}
\frac{(x^{p}+1)^{p^{k-1}\delta}-1}{p^k}
\equiv \sum_{\delta \in \Delta} \frac{(y+1)^{p^{k-1}\delta}-1}{p^k} \mod {\m_n}.$$ 
Hence we have $\sum_{\delta}\frac{(x+1)^{p^{k}\delta}-1}{p^k} \in \m_n + k_{n-1}$. 
Since $\ell(x)$ is convergent, for sufficiently large $k_0$, we have
$\sum_{k=k_0}^{\infty}\sum_\delta \,
\frac{(x+1)^{p^{k}\delta}-1}{p^k} \in \m_n$, and  therefore 
$\ell(x)$ is contained  in $ \m_n + k_{n-1}$. 
Since $\ell$ is injective on $\mathcal{F}(\m_n)(\cong \widehat{\G}_m(\mathfrak{m}_n) )$ 
 and 
is compatible with the Galois action,  
we have $\ell \left( \m_n\right)  \cap k_{n-1}=\ell \left( \m_{n-1}\right) $. 
Therefore we have an injection 
$$\ell (\m_n)  /\ell (\m_{n-1}) \hookrightarrow 
(\m_n+k_{n-1})/k_{n-1} \cong \m_n/\m_{n-1}.$$ 
By direct calculations, we have
$\ell (\iota^{-1}(c_n)) \equiv  \sum_{\delta} (\zetapnpi^\delta-1) \mod 
{k_{n-1}}$. Since $ \sum_{\delta} (\zetapnpi^\delta-1)$ generates 
$\m_n/\m_{n-1}$ as a $\Z_p[\Gamma_n]$-module with respect to the usual addition, 
the above injection is in fact a bijection.  
Thus $(\iota^{-1}(c_n)^\sigma)_{\sigma \in \Gamma_n }$  generate 
$\mathcal{F}(\mathfrak{m}_n)/\mathcal{F}(\mathfrak{m}_{n-1})$ and 
by induction $(\iota^{-1}(c_n)^\sigma)_{\sigma \in \Gamma_n }$ and $\varepsilon$ generate $\mathcal{F}(\mathfrak{m}_n)$. 
Since $\widehat{\G}_m$ is isomorphic  to $\mathcal{F}$ by $\iota$, 
we have ii). 
\end{proof}

Since $Nd_n=d_0=1$, by Hilbert's theorem 90, there exists an element $x_n \in k_n$ such that 
$d_n=x_n^{\gamma}/x_n$ for a fixed generator $\gamma$ of $\Gamma$.  
We put $\pi_n=\prod_{\delta \in \Delta} (\zetapnpi^\delta-1)$. Then 
$\pi_n$ is a uniformizer of $k_n$ and by the previous proposition, $x_n$ can be taken 
of the form $x_n=\pi_n^{e_n} u_n$ where $u_n \in (U^1_n)^{\mathrm{N}=1}$. 

\begin{proposition} \label{key prop}
In the same notation as the above, we have 
$$p \,\equiv\, e_n \,(p-1)\, \log_p\kappa(\gamma) \mod p^{n+1}.$$
\end{proposition}
\begin{proof}
If we put $$G(X)=\exp(p) \cdot \exp \circ\; \ell\;(X) = \exp \circ\; \ell\;(X[+] \varepsilon)\in 1+(p, X)\Zp[[X]],$$ then  by definition 
$$G_\sigma(\zetapmpi-1)=d_m^{\sigma}$$ where 
$G_\sigma(X)=G((X+1)^{\kappa(\sigma)}-1)$ for $\sigma \in \Gamma$. 
By Proposition \ref{norm} ii), $u_n$ is written as a combination of $d_n^\sigma$, and 
using the above $G_\sigma$, we can construct  a power series $H(X)  \in 1+(p, X)\Zp[[X]]$
such that $H(\zetapmpi-1)=N_{k_n/k_{m}} u_n$ for $0 \leq m \leq n$. 
We put 
$$F(X)=\left(
\prod_{\delta \in \Delta} \frac{(X+1)^{\delta \kappa(\gamma)}-1}{(X+1)^\delta-1}\right)^{e_n}\frac{H((X+1)^{\kappa(\gamma)}-1)}{H(X)}.$$
Then we have 
$$G(X) \equiv F(X)\quad  \mod \frac{(X+1)^{p^{n+1}}-1}{X}$$ 
since they are equal if we substitute $X=\zetapmpi-1$ for $0 \leq m \leq n$. 
Substituting $X=0$ in this congruence and 
taking the $p$-adic logarithm, we see that  $p \equiv e_n (p-1) \log_p\kappa(\gamma) \mod p^{n+1}.$
\end{proof}

\section{The Coleman map for the Tate curve.}

We construct the Coleman map for the Tate curve 
following the Appendix of  \cite{R}  or Section 8 
of \cite{K}. See also \cite{Ku}. 
In this section we assume that $E$ is the Tate curve  
$$E_q : y^2+xy=x^3+a_4(q)x+a_6(q)$$ 
where $q=q_E \in \Q_p^\times$ satisfying $|q|_p<1$ and 
$$s_k(q)=\sum_{n \geq 1} \frac{n^kq^n}{1-q^n}, \quad a_4(q)=-s_3(q), \quad a_6(q)=-\frac{5s_3(q)+7s_5(q)}{12}.$$ 
Then we have the uniformization 
$$\phi:\; {\C_p}^\times/q^\Z \cong E_q({\C_p}), 
\qquad u \mapsto (X(u,q), Y(u,q))$$ 
where 
$$X(u,q)=\sum_{n \in \Z} \frac{q^nu}{(1-q^nu)^2}-2s_1(q),$$
$$Y(u,q)=\sum_{n \in \Z} \frac{(q^nu)^2}{(1-q^nu)^3}+s_1(q).$$
(Of course, we put $\phi(q^\Z)=O$.) This isomorphism induces the 
isomorphism of the formal groups  $\widehat{\phi} : \widehat{\G}_m \cong \widehat{E}$. 
It is straightforward to see that the pull back of the invariant differential 
$\omega_E=\frac{dx}{2y+x}$  on $\widehat{E}$ with the parameter 
$t=-x/y$ by $\widehat{\phi}$ is the 
invariant differential 
$\omega_{\widehat{\G}_m}=\frac{dX}{1+X}$ on $\widehat{\G}_m$ 
with the parameter $X=u-1$. 
Hence  
$\widehat{\phi}$ is given by the power series $t=\exp_{\widehat{E}} \circ \log (1+X)-1 
\in \Z_p[[X]]$.  
%and its inverse  is given by $\exp \circ \log_{\widehat{E}} (X)-1$.  

From now we identify  $\widehat{\G}_m $ with $\widehat{E}$ by $\widehat{\phi}$. 
In particular, we regard $c_n \in \widehat{\G}_m (\mathfrak{m}_n)$ in the previous section
 as an element of  $\widehat{E}(\mathfrak{m}_n)$.  

Let $T=T_pE$ be the $p$-adic Tate module of $E$ and 
$V=T\otimes \Q_p$. 
The cup product  induces a non-degenerate pairing of 
Galois cohomology groups 
\begin{equation*}\label{sec: const, eq: pairing}
(\;,\;)_{E,n}: \; H^1(k_n, T)\times H^1(k_n, T^*(1)) 
\rightarrow H^2(k_n, \Zp(1)) \cong \Zp.
\end{equation*}
If there is no fear of confusion, we write $(\;,\;)_{E,n}$ simply as $(\;,\;)_{E}$.  
By the Kummer map, we regard $\widehat{E} (\m_n)$ as a subgroup of  $H^1(k_n, T)$.  
Then we define a morphism 
$\mathrm{Col}_{n}: H^1(k_n, T^*(1)) \rightarrow \Zp[\Gamma_n]$ by 
$$ z \quad \longmapsto 
\sum_{\sigma \in \Gamma_n} \;( c_n^{\sigma}, \;z)_{E,n} \; \sigma. $$  
This morphism is compatible with the natural Galois action and 
since the sequence $(c_n)_n$ is norm compatible, $\mathrm{Col}_n$ is also compatible for $n$ 
with respect to the corestrictions and the natural projections. 
We define the Coleman map 
$$\mathrm{Col}: \quad \varprojlim_n H^1(k_n, T^*(1)) \;\longrightarrow\; \Lambda=\Z_p[[\Gamma]]$$
as the projective limit of $\mathrm{Col}_n$ over all $n$. 

We recall the dual exponential map. 
For every $n$ let $\mathrm{tan}(E/k_n)$ denote the 
tangent space of $E/k_n$ at the origin, and consider the Lie group 
exponential map 
$$ \exp_{E,n} : \;\mathrm{tan}(E/k_n) \rightarrow E(k_n) \otimes \Q_p.$$ 
The cotangent space 
$\mathrm{cotan}(E/k_n)$ is generated by the 
invariant differential $\omega_E$ over $k_n$, and 
we let $\omega_E^*$ be the corresponding dual basis of $\mathrm{tan}(E/k_n)$. 
Then there is a dual exponential map
$$\exp_{E,n}^*:\;
 H^1(k_n, V^*(1)) \longrightarrow \mathrm{cotan}(E/k_n)=k_n\,\omega_E,$$
which has a property 
\begin{equation*}\label{sec: const, eq: dualexp}
(x, z)_{E,n} = \mathrm{Tr}_{k_n/\Q_p} \log_{\widehat{E}} (x) \,\exp_{\omega_E,n}^*(z)
\end{equation*}
for every $x \in \hat{E}(\m_n)$ and $z \in H^1(k_n, V^*(1))$. 
Here 
$\exp_{\omega_E,n}^*=\omega_E^* \circ \exp_{E,n}^*$. 
If there is no fear of confusion, we write $\exp_{\omega_E,n}^*(z)$ as $\exp_{\omega_E}^*(z)$. 
Then using the identification $\widehat{\phi} : \widehat{\G}_m \cong \widehat{E}$,
 the morphism $\mathrm{Col}_n$ is described 
in terms of the dual exponential map as follows. 

\begin{align*}
\mathrm{Col}_{n} (z)
&=\sum_{\sigma \in \Gamma_n} \;( c_n^{\sigma}, z)_{E,n} \; \sigma \\
&=\sum_{\sigma \in \Gamma_n} 
(\,\mathrm{Tr}_{k_n/\Q_p} \log_p(d_n^\sigma)\, 
\exp_{\omega_E}^*(z)\,)\,\sigma \\
&=\left(\sum_{\sigma \in \Gamma_n}\,\log_p(d_n^\sigma)\,\sigma\right)\;
\left(\sum_{\sigma \in \Gamma_n}\,\exp_{\omega_E}^*(z^\sigma)\,\sigma^{-1}\right).
\end{align*}

Let $G_n$ be the Galois group $\mathrm{Gal}(\Q_p(\zetapn)/\Q_p)$ and 
let $\chi$ be a finite character of $G_{n+1}$ of conductor $p^{n+1}$ which is trivial on  $\Delta$. 
Then we have 
\begin{align*}
\sum_{\sigma \in \Gamma_n}\,\log_p(d^\sigma_{n})\,\chi(\sigma)=
\begin{cases}
\tau(\chi)
&\:\: \text {if $\chi$ is non-trivial}, \\
\quad 0 &\:\: 
\text {otherwise}\\
\end{cases}
\end{align*}
where $\tau(\chi)$ is the Gauss sum
$\sum_{\sigma \in G_{n+1}}\chi(\sigma)\, \zeta_{p^{n+1}}^{\sigma}$.  
Hence for $\chi \not=1$, we have 
$$\chi \circ \mathrm{Col}(z)=\tau(\chi) \sum_{\sigma \in \Gamma_n}\,\exp_{\omega_E}^*(z^\sigma)\,\chi(\sigma)^{-1}.$$
Kato showed that there exists an element 
$z^{\mathrm{Kato}}   \in \varprojlim_n H^1(k_n, T^*(1))$ such that 
\begin{align*}
\sum_{\sigma \in \Gamma_n}\,\exp_{\omega_E}^*((z^\mathrm{Kato})^\sigma)\,\chi(\sigma)^{-1}=
%\begin{cases}
 e_p(\overline{\chi})\frac{L(E, \overline{\chi}, 1)}{\Omega^+_E}
%&\:\: \text {if $\chi$ is non-trivial}, \\
% \left(1-\frac{1}{p}\right)  \frac{L(E, 1)}{\Omega^+_E}&\:\: 
%\text {otherwise}. \\
%\end{cases}
\end{align*}
where $e_p(\chi)$ is the value at $s=1$ of the $p$-Euler factor of $L(E, {\chi}, s)$, 
that is, $e_p(\chi)=1$ if $\chi$ is non-trivial and $e_p(\chi)=\left(1-\frac{1}{p}\right) $ if 
$\chi$ is trivial. 
(See \cite{Ka}, Theorem 12.5.) Hence we have 
$$ 
\chi \circ \mathrm{Col}(z^{\mathrm{Kato}})=\tau(\chi) \frac{L(E, \overline{\chi}, 1)}{\Omega^+_E}$$ 
if $\chi$ is non-trivial. The $p$-adic $L$-function $L_p(E,s)$ is written of the form 
$$L_p(E,s)=\mathcal{L}_{p, \gamma}(E, \kappa(\gamma)^{s-1}-1)$$
for some power series $\mathcal{L}_{p, \gamma}(E, X) \in \Z_p[[X]]$. 
If we identify $\Lambda=\Zp[[\Gamma]]$ with $\Zp[[X]]$ by sending $\gamma \mapsto 1+X$,  
then it satisfies  an  interpolation formula   
$$ 
\chi \circ\mathcal{L}_{p, \gamma}(E,X) =\tau(\chi) \frac{L(E, \overline{\chi}, 1)}{\Omega^+_E}.$$
Since an element of $\Lambda$ has only finitely many zeros, 
we conclude that  
$$\mathrm{Col}(z^{\mathrm{Kato}})(X)=\mathcal{L}_{p, \gamma}(E,X).$$
Here we denote  $\mathrm{Col}(z^{\mathrm{Kato}})$ by 
$\mathrm{Col}(z^{\mathrm{Kato}})(X)$ to emphasis 
that we regard $\mathrm{Col}(z^{\mathrm{kato}})$ as a power series in $\Z_p[[X]]$. 
Note that  we have 
$\bold{1} \circ \mathrm{Col}(z)=0$  for the trivial character $\bold{1}$, or
$\mathrm{Col}(z)(0)=0$,  namely, any Coleman power series $\mathrm{Col}(z)(X)$ for 
the Tate curve 
 has a trivial zero at $X=0$. 

\section{The first derivative of the Coleman map.} 

We compute the first derivative of the Coleman map $\mathrm{Col}(z)(X)$. 
By Tate's uniformization,  
there is an exact sequence of local Galois representations 
\begin{equation}\label{exact}
0 \rightarrow T_1 \rightarrow T \rightarrow T_2 \rightarrow 0
\end{equation}
where $T_1=T_p\widehat{E} \cong \Zp(1)$ and  $T_2 \cong \Z_p$. 
The cup product 
induces  a non-degenerate paring 
\begin{equation*}\label{sec: const, eq: pairing}
 H^1(k_n, T_1)\times H^1(k_n, T_1^*(1)) 
\rightarrow H^2(k_n, \Zp(1)) \cong \Zp.
\end{equation*}
With the identification by  $\widehat{\phi}: T_1 \cong  \Zp(1)$, 
this is in fact the cup product 
pairing of $\mathbb{G}_m$
\begin{equation*}\label{sec: const, eq: pairing}
(\;,\;)_{\mathbb{G}_m,n}: \; H^1(k_n, \Z_p(1))\times H^1(k_n, \Zp) 
\rightarrow H^2(k_n, \Zp(1)) \cong \Zp.
\end{equation*}
If there is no fear of confusion, we write $(\;,\;)_{\mathbb{G}_m,n}$ simply as  $(\;,\;)_{\mathbb{G}_m}$. 
Since $c_n \in \widehat{E} (k_n) \subset  H^1(k_n, T_1)$, 
 we have 
$$(c_n^\sigma,\; z)_{E,n}=(d_n^\sigma, \,\pi(z))_{\mathbb{G}_m, n}$$
 for $z \in H^1(k_n, T^*(1))$ 
where  $\pi$ is the morphism induced by the projection $T^*(1)  \rightarrow T_1^*(1)$. 
Tate's uniformization $\phi$ also induces a commutative diagram  
\begin{equation*}
\begin{CD}
 H^1(k_n, V^*(1)) @>\exp_E^*>> k_n\,\omega_E @>\omega_E^*>> k_n\\
@V \pi VV @. @VVV\\
 H^1(k_n, V_1^*(1)) @>\exp_{\mathbb{G}_m}^*>> k_n\,\omega_{\mathbb{G}_m}
@>\omega_{\mathbb{G}_m}^*>> k_n
\end{CD}
\end{equation*}
where $\omega_{\mathbb{G}_m}$ is the invariant differential  of $\mathbb{G}_m$
which is $\frac{dX}{1+X}$ on  $\widehat{\mathbb{G}}_m$, and 
$\omega_{\mathbb{G}_m}^*$ is the dual basis for $\omega_{\mathbb{G}_m}$.  
We also put $\exp^*_{\omega_{\mathbb{G}_m}}
=\omega_{\mathbb{G}_m}^* \circ\,\exp^*_{\omega_{\mathbb{G}_m}}$. 

Now we compute the derivative. With the same notation as the previous section, 
we have 
\begin{align*}
\mathrm{Col}_n(z) &=\sum_{\sigma \in \Gamma_n} \;(c_n^{\sigma},\; z)_{E,n} \; \sigma =\sum_{\sigma \in \Gamma_n} \;(d_n^{\sigma}, \,\pi(z))_{\mathbb{G}_m,n} \; \sigma \\
&=\sum_{\sigma \in \Gamma_n} \;(( x_n^{\gamma}/x_n)^{\sigma},\, \pi(z))_{\mathbb{G}_m,n} \; \sigma \\
&=(\gamma^{-1}-1)\sum_{\sigma \in \Gamma_n} \;( x_n^{\sigma},\, \pi(z))_{\mathbb{G}_m,n} \; \sigma. 
\end{align*}
Therefore by the identification $\Zp[X]/((X+1)^{p^{n}}-1) \cong \Zp[\Gamma_n]$, 
$X \mapsto \gamma-1$, we have 
$$\frac{\mathrm{Col}(z)(X)}{X} \;\equiv\; -\frac{1}{\gamma } \sum_{\sigma \in \Gamma_n} \;(x_n^{\sigma}, \,\pi(z))_{\mathbb{G}_m,n} \; \sigma \mod \frac{(X+1)^{p^{n}}-1}{X}. 
$$
Hence 
$$\mathrm{Col}(z)'(0)\;\equiv\; - \;(\mathrm{N} x_n, \;\pi(z))_{\mathbb{G}_m,0} \;  \mod p^n. 
$$
Since $\mathrm{N} x_n=p^{e_n} \mathrm{N}(u_n)=p^{e_n}$ and by Proposition \ref{key prop},  we have 
{
$$(\mathrm{N} x_n, \,\pi(z))_{\mathbb{G}_m} = e_n\, (p,\,\pi(z))_{\mathbb{G}_m} 
\equiv \frac{p}{(p-1)\log_p \kappa(\gamma)} \;(p,\, \pi(z))_{\mathbb{G}_m}  \mod p^n.$$ }
Taking limit for $n$, we see that 
\begin{equation}\label{key}
\mathrm{Col}(z)'(0)=-\frac{p}{(p-1)\log_p \kappa(\gamma)}\;(p,\,\pi(z))_{\mathbb{G}_m}.  
\end{equation}

Next we compute $(p,\,\pi(z))_{\mathbb{G}_m}$. 
We consider the exact sequence 
\begin{equation*}
\begin{CD}
H^1(\Q_p, T^*(1)) @>\pi>>H^1(\Q_p, T^*_1(1)) @>\delta_2>> H^2(\Q_p, T_2^*(1))
\end{CD}
\end{equation*} induced by (\ref{exact}), 
and a  diagram 
\begin{equation*}
\begin{CD}
H^1(\Q_p, T_1) @. \;\times\; H^1(\Q_p, T^*_1(1)) @>(\,,\,)_{\mathbb{G}_m}>> H^2(\Q_p, \Zp(1))=\Zp \\
@A\delta_1AA  @V\delta_2VV @VVV \\
H^0(\Q_p, T_2) @. \;\times\; H^2(\Q_p, T^*_2(1))  @>(\,,\,)_{\mathbb{G}_m}>> H^2(\Q_p, \Zp(1))=\Zp.
\end{CD}
\end{equation*}

It is straightforward to see that the connecting morphism $\delta_1$ is given by 
$$H^0(\Q_p, T_2)=\Zp \;\rightarrow\; \Q_p^\times \otimes \Zp=H^1(\Q_p, T_1) , \quad 
1 \,\mapsto\, q_E\otimes 1.$$ 
Hence for $w \in H^1(\Q_p, T^*_1(1))$, we have 
$$(q_E\otimes 1,\,w)_{\mathbb{G}_m}=(\delta_1(1),\,w)_{\mathbb{G}_m}
=(1,\,\delta_2(w))_{\mathbb{G}_m}.$$ 
In particular, if $w$ comes from $H^1(\Q_p, T^*(1))$, 
namely,  it is of the form $\pi(z)$,  then 
\begin{equation}\label{q=0}
(q_E\otimes 1,\,w)_{\mathbb{G}_m}=(q_E\otimes 1,\,\pi(z))_{\mathbb{G}_m}=0. 
\end{equation}
On the other hand, if we put $q_E=p^{\mathrm{ord}_p (q_E)}\,\rho\, u_{q}$ where $ \rho \in \mu_{p-1}$ and 
$u_{q} \in 1+p\Zp$, we have 
\begin{align}
(q_E\otimes 1,\,w)_{\mathbb{G}_m}&=\mathrm{ord}_p (q_E)\;(p,\,w)_{\mathbb{G}_m}
+(u_{q},\,w)_{\mathbb{G}_m}\\
& =\mathrm{ord}_p (q_E)\;(p,\,w)_{\mathbb{G}_m} \label{qcomp}
+\log_p (u_q) \, \exp^*_{\omega_{\mathbb{G}_m}}(w). 
\end{align}
Hence by (\ref{q=0}) and (\ref{qcomp}) we have 
\begin{equation}\label{key2}
(p,\,\pi(z))_{\mathbb{G}_m}=-\frac{\log_p (u_q)}{\mathrm{ord}_p (q_E)}
 \exp^*_{\omega_{\mathbb{G}_m}}(\pi (z))=-\frac{\log_p (q_E)}{\mathrm{ord}_p (q_E)}
 \exp^*_{\omega_{E}}(z).  
\end{equation}
Combining (\ref{key}) and (\ref{key2}),   we obtain 

\begin{theorem}\label{dcol}
For $z \in \varprojlim_n H^1(k_n, T^*(1))$, the leading coefficient of the Coleman map 
$\mathrm{Col}(z)$ is given by 
$$
\frac{d}{dX} \mathrm{Col}(z)(X)\;\vert_{{X=0}}=\frac{p}{(p-1)\log_p \kappa(\gamma)}\;
\frac{\log_p (q_E)}{\mathrm{ord}_p (q_E)}
 \exp^*_{\omega_{E}}(z).  $$
\end{theorem} 

Now  if $E/\Q$ has  split multiplicative reduction at $p$, then 
we may assume that $E$ is locally the Tate curve for some $q_E \in \Q_p^\times$.
We apply the above formula to Kato's element $z=z^{\mathrm{Kato}}$. 
Since   
$\exp^*_{\omega_{E}}(z^{\mathrm{Kato}})=(1-\frac{1}{p}) \frac{L(E,1)}{\Omega^+_E}$,  
we have 

\begin{corollary} 
Let $\mathcal{L}_{p,\gamma}(E, X)$ be the power series in $\Zp[[X]]$ such that 
$L_p(E,s)=\mathcal{L}_{p, \gamma}(E, \kappa(\gamma)^{s-1}-1)$. Then 
$$
\frac{d}{dX} \mathcal{L}_{p,\gamma}(E, X)\;\vert_{{X=0}}=\frac{1}{\log_p \kappa(\gamma)}\;
\,\frac{\log_p (q_E)}{\mathrm{ord}_p (q_E)}
\,\frac{L(E,1)}{\Omega^+_E}, 
$$
or 
$$
\frac{d}{ds} {L}_{p}(E, s)\;\vert_{{s=1}}=
\,\frac{\log_p (q_E)}{\mathrm{ord}_p (q_E)}
\,\frac{L(E,1)}{\Omega^+_E}.  
$$

\end{corollary}

%%--------------------Here the manuscript ends--------------------------------
\Addresses
\end{document}